\documentclass[12p, leqno]{article}
\usepackage{amsmath, amsthm}
\usepackage{amsfonts} 
\usepackage{amssymb}
\usepackage{graphicx}
\usepackage{color}

\newtheorem{thm}{Theorem}[section]

\newtheorem{prop}[thm]{Proposition}

\numberwithin{equation}{section}

\newcommand{\RR}{\mathbb{R}}
\newcommand{\ren}{\mathbb{R}^n}

\newcommand{\dy}{\,{\rm d}y}

\def\quotient#1#2{%
    \raise1ex\hbox{$#1$}\Big/\lower1ex\hbox{$#2$}%
}

\def\qed{\,\unskip\kern 6pt \penalty 500
\raise -2pt\hbox{\vrule \vbox to8pt{\hrule width 6pt
\vfill\hrule}\vrule}\par}

\definecolor{darkblue}{rgb}{0.05, .05, .65}
\definecolor{darkgreen}{rgb}{0.05, .70, .05}
\definecolor{darkred}{rgb}{0.8,0,0}
\begin{document}
\title{ \huge \bf Asymptotic behaviour  for the \\ Fractional Heat Equation \\ in the Euclidean space}
\author{\Large Juan Luis V\'azquez \\[8pt] Universidad Aut\'{o}noma de Madrid, Spain. } 
\date{} 

\maketitle

\begin{abstract}
We consider weak solutions of the fractional heat equation posed in the whole $n$-dimensional space, and establish their asymptotic  convergence to the fundamental solution as $t\to\infty$ under the assumption that the initial datum is an integrable function, or a finite Radon measure. Convergence with suitable rates is obtained for solutions with a finite first initial moment, while for solutions with compactly supported initial data convergence in relative error holds. The results are applied to the fractional Fokker-Planck equation. Brief mention of other techniques and related equations is made.
\end{abstract}

\vspace{1cm}


\normalcolor

\normalsize
\section{Introduction}

We consider weak solutions of the fractional heat equation (FHE)
\begin{equation}\label{FHE}
\partial_t u + (-\Delta)^s u=0, \quad 0<s<1\,,
\end{equation}
posed in the  $n$-dimensional space $\ren$, $n\ge 1$.
The fractional Laplace operator $(-\Delta)^s$ may be defined through its Fourier transform, or by its representation
$$
(-\Delta)^s f(x)=c(n,s)\int_{\ren} \frac{f(x)-f(y)}{|x-y|^{n+2s}}\,\dy
$$
with $0<s<1$, see its basic properties in \cite{Landkoff, Stein}. In the limit $s\to 1$ the standard Laplace operator, $-\Delta$, is recovered  (see Section 4 of \cite{DNPV2012}), but there is a big difference between the local operator $-\Delta$ that appears in the classical heat equation and represents Brownian motion, and the nonlocal family $(-\Delta)^s $, $0<s<1$. The latter are generators of  L\'evy processes that include jumps and long-distance interactions, resulting in anomalous diffusion \cite{App, Bert}.

Applying Fourier transform to the FHE (with respect to the space variable $x$, and calling the new variable $\xi$), we get the equation \ $\partial_t \widehat{u}=-|\xi|^{2s}\widehat{u}$,
that allows to solve the initial-value problem in Fourier space by means of the formula
\begin{equation}
\widehat{u}(\xi,t)=\widehat{u_0}(\xi)\,e^{-|\xi|^{2s}t}\,.
\end{equation}
Applying the inverse transform, the fractional heat equation can be solved for all $0<s\le 1$ by means of the fundamental solution, $P(x,t;s)$ which is the inverse transform of  the function $f_s(\xi)=e^{-|\xi|^{2s}t}. $ We have the representation formula
\begin{equation}
u(x,t)=\int_{\ren} P(x-y,t;s)u_0(y)\,dy.
\end{equation}
The optimal class of initial data in this representation is given  by Bonforte, Sire and the author in  \cite{BSV2017}, and is the class of locally finite Radon measures satisfying the growth condition
\begin{equation}
\int_{\ren}  (1+|x|^2)^{-(n+2s)/2}\,d|\mu|(x)<\infty.
\end{equation}
In the present paper the initial data will integrable $d\mu(x)=u_0(x)dx$ with $u_0\in L^1(\ren)$, and we will use this formula to obtain the asymptotic behaviour of solutions with finite space integral, otherwise known as finite mass solutions.

There are a number of previous works dealing with the FHE, in particular in the probabilistic theory, cf. \cite{BG1960, BJ2007, ChenKum2008, MK2000, Vald}.  From the analytical point of view we may mention paper \cite{BSV2017} that contains a rather general theory of existence, uniqueness, initial traces, as well as a priori estimates and regularity of solutions, let us also mention the previous work \cite{Barrios2014}. Both papers contain references to the literature. For a comparative presentation of convergence to the  Gaussian kernel in the case of the classical heat equation, we refer to \cite{Vaz17}, as well as many classical references.

\medskip

\noindent {\bf Outline of the contents.} After a  section devoted to establish the needed properties of the heat kernel, the main result about convergence to the fundamental solution is carefully stated as Theorem \ref{main.convthm} in Section \ref{sec.conv}. Also the stronger convergence under the assumption of finite first initial moment is stated, Theorem \ref{thm.conv.2}. Proofs of both results are given in Subsection \ref{sec.conv1}. A detailed list of comments on the convergence results is examined in Subsection \ref{sec.comments}. Let us select one: even if finer convergence is obtained under the assumption of finite first initial moment, the first moment is never finite for solutions with compactly supported data at any positive time if $0<s\le 1/2$.

Section \ref{conv.re} discusses the topic of convergence in relative error, a quite strong result that is not true in classical heat equation, $s=1$ (the Gaussian case).

Section \ref{sec.ffpe} deals with the application of our results to the Fractional Fokker-Planck Equation, which is a quite interesting and active topic. A brief mention of the Fractional Ornstein-Uhlenbeck equation is made.

Finally, Section   \ref{sec.other} contains some ideas on convergence for related nonlinear equations that have been recently investigated or are currently being investigated.


\section{Properties of the fractional heat kernel}
\label{sec.prop.hk}

Before proceeding with the main topic of the paper, it will be convenient to study in some detail the properties of the fundamental solution since it serves as kernel of the representation formula.

\medskip

\noindent $\bullet$ The kernel is explicit for two particular values, $s=1$ and $s=1/2$. In the first case we get the Gaussian kernel for the standard heat equation:
\begin{equation}
P(x,t;1)=(4\pi t)^{-n/2}\,e^{-|x|^2/4t}\,,
\end{equation}
which is selfsimilar with $P=t^{-n/2}F_1(xt^{-1/2})$ with profile $F_{1}(x)=(4\pi)^{-n/2}\,e^{-|x|^2/4}$. Notice the very fast decay with the typical square exponential tails as $|x|\to\infty$.
In the fractional case $s=1/2$ we have
\begin{equation}
 P(x,t;1/2)=\frac{c_n t}{(t^2+| x|^2)^{(n+1)/2}}\,, \quad \mbox{with profile} \ F_{1/2}(x)=\frac{c_{n}}{(1+| x|^2)^{(n+1)/2}},
\end{equation}
where $c_n=\Gamma(n+(1/2))/\pi^{n+(1/2)}$ is a normalization constant. The decay in $|x|$  is now slower with a decay power that depends directly on $s$. It is interesting to point out that this kernel is precisely the Poisson kernel for the Laplace equation posed in the upper half-space
$H^+=\{(x,y): x\in \ren, y>0\}\subset \RR^{n+1}$, after replacing the last coordinate $y$ by time $t$. This fact is easily justified using the Caffarelli-Silvestre extension theorem \cite{CaffSilv}.  The example  $s=1/2$ will be the paradigm for the rest of the fractional cases $0<s<1$ as we will see below.

\medskip

\noindent $\bullet$  For general $0<s<1$ the fortunate properties of the case $s=1/2$ do not exist and do not have explicit formulas for the kernel. However, interesting information is known and it will be enough for our purposes. Thus, it is well-known that the kernel has the self-similar form
\begin{equation}
P(x,t;s)=t^{-n/2s}F_s(|x|t^{-1/2s})\,.
\end{equation}
Since the Fourier symbol $ e^{-t|\xi|^{2s} }$ is a tempered distribution, it immediately follows  that $P(x,t;s) \in C^\infty ((0,\infty)\times\ren )$, see  \cite{BJ2007, Kol2000}.
The asymptotic formula
$$
\lim_{|\xi|\to\infty}|\xi|^{n+2s}F_s(\xi) = C(n,s)\,,
$$
was established by Blumental and Getoor \cite{BG1960} in 1960 by a delicate analysis based on integrals in the complex plane, and it generalizes in a precise way the power-like tail behaviour observed in the case $s=1/2$. This means that for all $x\in\ren$ and  $t>0$
\begin{equation}\label{HK-estimates.1}
P(x,t;s) \asymp  \frac{t}{\big(t^{1/s}+|x|^2\big)^{(n+2s)/2}}\,,
\end{equation}
where $\asymp$ means that the ratio bounded by a constant factor from above and below.

\medskip

\noindent $\bullet$  Taking the selfsimilar form $P=t^{-n/2s}F(|x|t^{-1/2s})$ and plugging it into the equation we find that all these profiles satisfy:
$$
2st\,\partial_t P= -t^{-n/2s}(nF+ r\partial_{r} F), \quad nP+ r\partial_r P= t^{-n/2s}(nF+ r\partial_r F)\,,
$$
which in other words means that
\begin{equation}\label{eq.fundsol}
2s\,(-\Delta)^s F = nF+ r\partial_{r} F.
\end{equation}

\noindent $\bullet$  Monotonicity in $r$. The property $\partial_r P\le 0$  follows by an application of Aleksandrov's reflection principle. The Aleksandrov reflection method is a well-established tool to prove monotonicity of solutions of wide classes of (possibly nonlinear) elliptic and parabolic equations
\cite{Aleks, Ser71}. The use for nonlinear parabolic equations is carefully explained in Chapter 9 of \cite{Vazpme}.  Use of the method for fractional operators is recent. We have  explained the application to (possibly nonlinear) fractional parabolic equations in \cite{VazBar2014}, Section 15, and we are here in such a situation. One of simplest consequences of the principle is that uniquely defined  fundamental solutions must be monotone decreasing in the radial space variable for every time $t>0$.

\medskip

\noindent $\bullet$  Boundedness of $\partial_t P$. The bound from below $2stP_t\ge -nP$ holds by the previous formulas and the monotonicity of $P$ in $r$. This means that
$$
\int (\partial_t P)^{-}dx\le \frac{n}{2st}\int P dx = \frac{n}{2st}\,,
$$
Using the fact that $\int P_t\,dx=0$ we get
$$
\int |\partial_t P|dx\le  \frac{n}{st}.
$$
Since $\partial_t P$ satisfies the FHE, is has an a priori estimate in $L^\infty$ (see  \cite{BSV2017}) and
$$
|\partial_t P(x,1)|\le  K(n,s)\,,
$$
where the constant is related to the constant of the bound for $P(x,t;s)$.

\medskip

\noindent $\bullet$  A finer analysis using the Fourier transform representation allows us to  obtain the asymptotic behaviour of $ \partial_t P=-(-\Delta)^{s} P_t$.   We can use Lemma 1, page 304, of Pruitt-Taylor \cite{Pruitt-Taylor} that  implies in particular that
$$
|(-\Delta)^{s} P(x,1;s)|\le C (1+|x|^2)^{-(n+2s)/2}\,.
$$
This is quoted and used in \cite{V4}. The estimate implies an interesting estimate for the time derivative,
\begin{equation}
t\,|\partial_t P(x,t;s)|/P(x,t;s) \le C_{n,s}\,.
\end{equation}
 The bounds for $\partial_t P(x,1)$ imply similar bounds for $r\partial_r F$. In particular
\begin{equation}
|r \partial_r F(r)|\le C(1+r^2)^{-(n+2s)/2}.
\end{equation}

All these formulas can be easily checked on the explicit example of case $s=1/2$, while they differ in a marked way from what happens for the classical heat equation, $s=1$. That fact makes for a common theory in all fractional cases  $0<s<1$.

%
\section{Asymptotic convergence to the Fundamental Solution}\label{sec.conv}

The main result on the asymptotic behaviour of general integrable solutions of the fractional heat equation (FHE) consists in proving that they look increasingly like the fundamental solution $P(x,t;s)$, here denoted as $P_t$ for the sake of brevity.  We start with a consideration that also applies in the Gaussian case: since $P_t$ goes to zero uniformly with time, the estimate of the convergence of $u$ to $P_t$ has to take into account that fact and compensate for it. This happens by considering a renormalized error that divides the standard error in some norm by the size of the fundamental solution $P_t$ in the same norm. For instance, in the case of the sup norm we know that
$$
\|P_t\|_{L^\infty(\ren)}=Ct^{-n/2s},\qquad \|P_t\|_{L^1(\ren)}=1.
$$
This is the basic result that we will prove

\begin{thm} \label{main.convthm} Let $u_0\in L^1(\ren)$ and let  $\int u_0(x)dx=M$ be its mass and let $M\ne 0$. Let $u(t)=u(\cdot,t)$ be the solution of the FHE with $s\in (0,1)$, posed in the whole space. We show that $u(t)$ ends up by looking like $M$ times the fundamental solution  $P_t=P(x,t;s)$ in the sense that
\begin{equation}\label{cr.l1}
\lim_{t\to\infty} \|u(t)-M P_t\|_1\to 0
\end{equation}
and also that
\begin{equation}\label{cr.linf}
\lim_{t\to\infty} t^{n/2s}\|u(t)-M P_t\|_\infty\to 0\,.
\end{equation}
By interpolation we get the convergence result for all $L^p$ norms
\begin{equation}\label{cr.lp}
\lim_{t\to\infty} t^{n(p-1)/2sp}\|u(t)-M P_t\|_{L^p(\ren)}\to 0\,.
\end{equation}
for all $1\le p\le \infty$. 
\end{thm}

Notice that we do not assume that $u$ is nonnegative, it may have any sign, and $M$ may be negative. In fact, the result also holds for $M=0$ but in that case we cannot say that $u(t)$ looks like a power of $P_t$.
About the three convergence results of the basic theorem, it is clear that \eqref{cr.lp} follows from \eqref{cr.l1} and \eqref{cr.linf}. Now, what is interesting is that \eqref{cr.linf} follows from \eqref{cr.l1} and the $L^1$ to $L^\infty$ effect (smoothing effect). We ask the reader to prove this fact. Hint: use \eqref{cr.l1} between $0$ and $t/2$ and then the  smoothing effect between $t/2$ and $t$. More comments will be given in Section \ref{sec.comments}.

The  proof we give here of such main result is based on the examination of the error in terms of the representation formula. This is very direct approach, but it needs a previous step,  whereby the proof  is done under the further restriction that the data have a finite first moment. We introduce the first moment of the distribution $u_0$ as the vector $\vec{\mathcal N}_1(u_0)$ with components
\begin{equation}
{\mathcal N}_{1,i}(u_0):=\int_{\ren} u_0(y)y_i\,dy\,, \qquad 1\le i\le n.
\end{equation}
 In the calculations we will need the first absolute moment
\begin{equation}
{\mathcal N}_1(u_0)=\int_{\ren} |u_0(y)y|\,dy <    \infty\,,
\end{equation}

Under this extra assumption the convergence result is more precise and quantitative. This case has an interest in itself since it shows the importance of controlling the first moment.
We already know that the first moment is associated to a conserved quantity.

\begin{thm} \label{thm.conv.2} Under the assumptions that $u_0\in L^1(\ren)$ and that the first absolute moment ${\mathcal N}_1(u_0)$, is finite,
we  get the convergence
\begin{equation}\label{conv.fm}
\displaystyle t^{n/2s}|u(x,t)-MP_t(x)|\le C_1 {\mathcal N}_1(u_0)t^{-1/2s}\,,
\end{equation}
as well as
\begin{equation}
\displaystyle \|u(x,t)-MP_t(x)\|_{L^1(\ren)}\le  C_2 {\mathcal N}_1(u_0)t^{-1/2s}\,.
\end{equation}
The rate $O(t^{-1/2s})$ is optimal under such assumptions.
\end{thm}

%
\subsection{Proof of the theorems via representation}\label{sec.conv1}

\noindent {\sc Proof of Theorem \ref{thm.conv.2}.}
(i)  We begin with a preliminary calculation
$$
\begin{array}{c}
\displaystyle u(x,t)-MP_t(x)= \int u_0(y)P_t(x-y)\,dy- P_t(x)\int u_0(y)\,dy\\[4pt]
\displaystyle =\int u_0(y)(P_t(x-y)-P_t(x))\,dy=
\int u_0(y)\left(\int_0^1\partial_\theta (P_t(x-\theta y))\,d\theta\right)\,dy\,.
\end{array}
$$
Using the selfsimilar formula for $P_t$ we get
$$
\displaystyle \displaystyle u(x,t)-MP_t(x)= - C t^{-n/2s}
\int dy \int_0^1 d\theta \,u_0(y) F'(|x-sy|/t^{1/2s})\langle y, \frac{x-\theta y}{|x-\theta y|}\rangle
t^{-1/2s}\,.
$$

(ii)  Let us do first the {\sl $L^1$ convergence.} We integrate in $x$ for fixed $t>0$ to get
$$
\begin{array}{c}
\displaystyle \|u(x,t)-MP_t(x)\|_1\le C \, t^{-N/2s}t^{-1/2s}
\int dx\int dy \int_0^1 d\theta \,|u_0(y) y| F'(|x-\theta y|/t^{1/2s}\\[4pt]
\displaystyle = C t^{-1/2s} \,\int dy \int_0^1 d\theta \,|u_0(y) y| \left(\int  t^{-n/2s}|F'(|x-sy|/t^{1/2s})| \,dx\right).
\end{array}
$$
With the change of variables $x-\theta y=t^{1/2s}\xi$ and putting $\xi=r\sigma$ in polar coordinates, we get
$$
\begin{array}{c}
\displaystyle \|u(x,t)-MP_t(x)\|_1\le C \, t^{-1/2s}
\iint dyd\theta \,|u_0(y) y| \int_{\ren}|F'(|\xi|)|d\xi\\[6pt]
\displaystyle =  C \, t^{-1/2s}
(\int_{\ren} |u_0(y) y|dy)(\int_0^1 d\theta) (\int_0^\infty |F'(r)|r^{n-1}dr)\,.
\end{array}
$$
Since $F'$ is smooth and has a sign we have
$$
\int_0^\infty |F'(r)|r^{n-1}dr= (n-1)\int_0^\infty F(r)r^{n-2}dr\,.
$$
In view of the known decay of $F$, see Section \ref{sec.prop.hk}, this integral is bounded, hence
$$
\|u(x,t)-MP_t(x)\|_1\le C \, t^{-1/2s} {\mathcal N}_1(u_0)\,.
$$

(iii) We show next the {\sl sup convergence}. Starting as before we get
$$
\displaystyle \displaystyle |u(x,t)-MP_t(x)| \le C t^{-(n+1)/2s}
\int dy \int_0^1 d\theta \,|u_0(y) y|\, |F'(|x-\theta y|/t^{1/2s})|\,.
$$
We know that $|F'|$ is bounded, hence
$$
\displaystyle |u(x,t)-MP_t(x)|\le  C_1t^{-(n+1)/2s}\int |u_0(y)y|\,dy\,.
$$
 Taking into account that $ P_t$ is of order $t^{-n/2}$ in sup norm, we write this result as \eqref{conv.fm}, and  $C=C(n,s)$ is a universal constant.

(iii) We could have  performed the proof under the further restriction that $u_0\ge 0$. For a signed solution we must then separate the positive and negative parts of the data and apply the results to both partial solutions. But we did not have to use the sign assumption.
\qed

\medskip

\noindent {\bf Optimality.} Take as solution the FS after a space displacement, $u(x,t)= P(x+he_1,t;s)$, and find the convergence rate to be exactly $O(t^{-1/2s})$. This is just a calculus exercise. {\sl Hint}: Write
$$
P(x+he_1,t)=P_t(x)+ h\partial_{x_1} P_t(\xi))
$$
(where $\xi=x+\theta h$, $0<\theta<1$), and check that
$$
|\partial_x P_t(\xi)| \le Ct^{-(n+1)/2s}|F'(|\xi|/t^{-1/2s}|=O(t^{-(N+1)/2s})\,,
$$
uniformly in $x$ small. \normalcolor This exercise shows that the term $h\partial_x P_t(x)$ is the {\sl precise corrector} with relative error $O(t^{-1/2s})$. We could continue the analysis by expanding in Taylor series with further terms.

\medskip

\noindent $\bullet$ {\sc Proof of Theorem \ref{main.convthm} using Theorem \ref{thm.conv.2}}. Given an initial function $u_0\in L^1(\ren)$ without any assumption on the first moment, we argue by approximation plus the triangular inequality. In the end we get a convergence result, but less precise. There will no specified rate of convergence.

(i) {\sl Proof of $L^1$ convergence.}  Let $u_0$ is integrable with integral $M$ and let us fix an error  $\delta>0$. First, we find an approximation  $u_{01}$ with compact support and such that
$$
\|u_0-u_{01}\|_1<\delta.
$$
Due to the  $L^1$ contraction property applied to the solution $u-u_1$, we know that for all  $t>0$
$$
\|u(t)-u_{1}(t)\|_1< \delta.
$$
On the other hand, we have just proved that for data with finite moment:
$$
\|u_1(t)-M_1P_t\|_1\le C {\mathcal N}_{1}(u_{0,\delta}) \,t^{-1/2s}\,.
$$
In this way, for sufficiently large $t$  (depending on $\delta$) we have
$$
\|u_1(t)-M_1P_t\|_1\le C \delta
$$
with $ C$ a universal constant. Next se recall that  $|M-M_1|\le \delta$. Using the triangular inequality, and letting then $\delta\to0$ we arrive at
\begin{equation}
\lim_{t\to\infty}\|u(t)-MP_t\|_\infty=0\,,
\end{equation}
which ends the proof.

(ii) {\sl Proof of $L^\infty$ convergence.} It follows from the previous one and the
the already mentioned effect $L^1\to L^\infty$,  applied to the solution $u-u_1$. Indeed, we know that for all  $t>0$
$$
\|u(t)-MP_t\|_\infty<C t^{-n/2s}\|u(t/2)-MP_{t/}\|_1\,.
$$

\noindent {\bf Counterexample.}  Let us explain how the lack of a rate of convergence for the whole class of integrable functions is shown. The idea is the same as in the Gaussian case treated in \cite{Vaz17}.

\begin{prop}Given any decreasing and positive rate function $\phi$ such that $\phi(t)\to 0$ as $t\to\infty$ we construct a modification of the Gaussian kernel that produces a solution with the same mass $M=1$ that satisfies the formula
$$
 t^{n/2s}\|u(x,t)- P_t(x)\|_\infty\ge k \phi(t_k)
$$
at a sequence of times $t_k\to\infty$ to be chosen.
\end{prop}

\medskip

\noindent {\sl Proof.}  Let us see the details. We have to find  a choice of small masses $m_1,m_2,\dots, $ with $\sum_k m_k=\delta< 1$, and locations $x_k$ with $|x_k|=r_k\to \infty$ and consider the solution
$$
u(x,t)=(1-\delta)P_t(x)+ \sum_{k=1}^\infty m_kP_t(x-x_k)\,.
$$
The error $u(x,t)-P_t(x)$  is calculated at $x=0$ as
$$
 t^{n/2s}|u(0,t)-P_t(0)|=|\delta F(0)-\sum_{k=1}^\infty m_k F(|x_k|/t^{1/2s})|=
\sum_{k=1}^\infty m_k (F(0)-F(|x_k|/t^{1/2s}))\,.
$$
We recall that for all large values of $r> 1$ we have
$$
F(r)\sim C(n,s)\,r^{n+2s}\,.
$$
Put $m_k=2^{-k}$ (any other summable series will do). Choose iteratively $t_k$ and $x_k$ as follows. Given choices for the steps $1,2,\dots, k-1$ , pick $t_k$ to be much larger than $t_{k-1}$ and such that $\phi(t_k)\le F(0) m_k/2k$. This is where we use that $\phi(t)$ tends to zero, even if it may decrease in a very slow way. Choose now $|x_k|=r_k$ so large  that $F(|x_k| t_k^{-1/2s}) < F(0)/2$. Essentially, the mass has to be displaced at distance equal or larger than $O(t_k^{1/2s})$. Then
$$
m_k(F(0)-F(|x_k|/t_k^{1/2s}))\ge m_k F(0) /2\ge k\phi(t_k)\,.
$$
This completes the proof. \qed


\subsection{Comments of the convergence results}
\label{sec.comments}

 We add important information to this result in a series of remarks.

\noindent $\bullet$  First, one comment about the spatial domain. The fact that we are working in the whole space is crucial for the result of Theorem \ref{main.convthm}. In the classical case $s=1$ the behaviour of the solutions of the heat equation posed in a bounded domain with different kinds of boundary conditions is also known and the asymptotic behaviour does not follow the Gaussian pattern. The same happens for $s<1$ with respect to the Fundamental Solution.

\noindent $\bullet$ The convergence to the Fundamental Solution happens on the condition that the data belong to the class of integrable functions. Well, this class can be extended a bit to the class of bounded Radon measures; this is actually no news since the theory \cite{BSV2017} says that the solution $u(t)$  corresponding to an initial measure $\mu\in {\mathcal M}(\ren)$ is integrable and bounded for any positive time, so we may change the origin of time and make the assumption of integrable and bounded data. But we point out the some of the proofs work directly for measures without any problem.

\noindent $\bullet$ It must be stressed that convergence to the FS {\sl  does not hold for other data}. Maybe  the simplest example of solution that does not approach the Gaussian is given by any non zero constant solution, but it could  be objected that $L^\infty(\ren)$ is very far from $L^1(\ren)$. Actually, the same happens for all $L^p(\ren)$ spaces, $p>1$. Indeed, a simple comparison theorem shows that for any $u_0\ge 0$ with $\int u_0(x)\,dx=+\infty$ we have
$$
\lim_{t\to\infty} t^{n/2s}u(x,t)=+\infty
$$
everywhere in $x\in \ren$ (and the divergence is locally uniform).

\noindent $\bullet$  The way different classes of non-integrable solutions actually behave for large time is an interesting question that we will not address here. Thus, the reader may prove using the convolution formula that
for locally integrable data that converge to a constant $C$ as $|x|\to\infty$ the solution $u(x,t)$ stabilizes  to that constant as $t\to\infty$. Taking growing data may produce solutions that tend to infinity with time. See details on growing solutions in \cite{BSV2017}, in particular classes of selfsimilar solutions are constructed.

\noindent $\bullet$ We recall that it is usually assumed that $u\ge 0$ on physical grounds but such assumption is not at all needed for the analytical study of this paper. As already mentioned, the basic result holds also for signed solutions even if the total integral is negative, $M\le 0$. There is no change in the proofs. We may also put $M=1$ by linearity as long as $M\ne 0$.

\noindent $\bullet$  Interesting case \  $M=0$: even if Theorem \ref{main.convthm} is true, the statement does not imply that the solution looks asymptotically like a Gaussian; to be more precise, it only says that the assumed first-order approximation disappears. If we want more precise details about what the solution looks like, we have to search further to identify the terms that may give us the size and shape of such a solution. This question will be addressed below. For the  moment let us point out that differentiation in $x_i$ of the Heat Kernel produces a new solution with zero integral
\begin{equation}
u_i(x,t)=\partial_{x_i} P_t(x)\,,
\end{equation}
where we can apply the above comments. In particular,
$$
t^{n/2s}u_i(t)= O(t^{-1/2s})\,,
$$
where $O(\cdot)$ is the Landau $O$-notation for orders of magnitude.

\noindent $\bullet$  There is a quite interesting comment to be made about Theorem \eqref{thm.conv.2}. It assumes that the initial function has finite first moment. This does not create any problems for $1/2<s<1$ since the fundamental solution has a finite first moment for all times. But when $0<s\le 1/2$ the absolute first moment of the fundamental solution is infinite and we wonder what kind of transition happens at $t=0+$. There a number of explanations to be given, but we leave this interesting question to the curiosity of the reader. Note the fundamental solutions have infinite second moment for all $0<s<1$. All this is in stark contrast to the classical case $s=1$ since all moments of the Gaussian function are finite.

\section{Convergence in relative error}\label{conv.re}

The factor $t^{n/2s}$ is the appropriate weight to consider relative error. We point out that the  relative error formula is given by
\begin{equation}
\epsilon_{rel}(t;u_0)=\frac{|u(x,t)-MP_t(x)|}{P_t(x)}
\end{equation}
An estimate on this quantity would be more precise than the previous theorems. It is known that the relative error does not admit a sup bound in the Gaussian case ($s=1$), as can be observed by choosing $u(x,t)=P_t(x-h)$ for some constant $h$ since then
$$
\epsilon_{rel}(t)=|e^{xh/2t}e^{-h^2/4t}-1|\,,
$$
which is not even bounded. This comment  shows that error calculations with Gaussians are delicate.

It is remarkable that we can do better in the fractional setting $0<s<1$ for all compactly supported data.

\begin{thm} \label{thm.conv.cs} Let us assume that $u_0\in L^1(\ren)$ and is compactly supported. Then we get the convergence with rate in relative error
\begin{equation}
 \epsilon_{rel}(t;u_0) \le CMR\,t^{-1/2s}
\end{equation}
for all large $t$, if the support of $u_0$ is contained in the ball of radius $R>0$.
\end{thm}

\noindent {\sl Proof.} We have
$$
\displaystyle \frac{ u(x,t)-MP_t(x)}{P_t(x)}= -\frac{C}{ F(|x|/t^{1/2s})}
\int dy \int_0^1 d\theta \,u_0(y) F'(|x-sy|/t^{1/2s})\langle y, \frac{x-\theta y}{|x-\theta y|}\rangle
t^{-1/2s}\,.
$$
Put $x=Xt^{1/2s}$, $y=Yt^{1/2s}$. and we put $B_R=B_R(0)$ we have
$$
\displaystyle \frac{ u(x,t)-MP_t(x)}{P_t(x)}= -Ct^{-1/2s}
\int_{ B_R}  \int_0^1  \,u_0(y) \frac{F'(|X-\theta Y|)}{ F(|X|)}\langle y, \frac{X-\theta Y}{|X-\theta Y|}\rangle dyd\theta\,.
$$
Hence,
$$
\left|\displaystyle \frac{ u(x,t)-MP_t(x)}{P_t(x)}\right|= CMt^{-1/2s} R
\sup_{\theta, Y}\frac{F'(|X-\theta Y|)}{ F(|X|)}\,.
$$
Since $0<\theta<1$ and $Y$ is bounded (even more, it goes to zero), and we know that
$F'(r)/F(r)\to 0$ as $r\to\infty$, the result follows. \qed

This result admits an improvement where the {\bf first-order correction} to the plain convergence result is also identified.

\begin{thm} \label{thm.conv.cs1} Let us assume that $u_0\in L^1(\ren)$ and is compactly supported. Then we get the improved convergence  in relative error with first corrector term:
\begin{equation}
 t^{1/2s}  \frac{|u(x,t)-MP_t(x) + {\mathcal N}_{1,i}(\partial_iP_t)(x)|}{P_t(x)}\to 0
\end{equation}
as $t\to\infty$, uniformly in $x\in \ren$. In other words, as $t\to\infty$
\begin{equation}
u(x,t)= M(u_0)\,P_t(x)-  {\mathcal N}_{1,i}(u_0)\,(\partial_iP_t)(x)+ o(t^{-1/2s})\,P_t(x)\,.
\end{equation}
\end{thm}

\noindent {\sl Proof.}  Let the support of $u_0$ is contained in the ball $B_R=B_R(0)$ of radius $R>0$. We have
$$
\begin{array}{c}
\displaystyle t^{(n+1)/2s}\,(u(x,t)-MP_t(x))) =
\displaystyle -C \int dy \int_0^1 d\theta \,u_0(y) F'(|x-sy|/t^{1/2s})\langle y, \frac{x-\theta y}{|x-\theta y|}\rangle \,,
\end{array}
$$
and
$$
\begin{array}{c}
\displaystyle t^{(n+1)/2s}\, \sum_i {\mathcal N}_i\partial_iP_t(x) =
\displaystyle C \int \,u_0(y) \,F'(x)\langle y, \frac{x}{|x|}\rangle dy\,.
\end{array}
$$
Combining them and putting  $x=Xt^{1/2s}$, $y=Yt^{1/2s}$, $e_i$ the $i$th component
of the unit vector $x/|x|$ and $e_i'$ the $i$-th component of $x- \theta y/|x- \theta y|$,  we get
$$
\begin{array}{c}
I(x,t):= \displaystyle \frac{ u(x,t)-MP_t(x) + \sum_i {\mathcal N}_i\partial_iP_t(x)}{P_t(x)}=\\ [10pt]
\displaystyle Ct^{-1/2s}  \sum_i
\int_{ B_R}  \int_0^1  \,u_0(y)y_i \,\frac{F'(|X|)e_i-F'(|X-\theta Y|)e'_i}{ F(|X|)}\,dyd\theta\,.
\end{array}
$$
Hence,
$$
t^{1/2s} \left|I(x,t)\right|= C  \sum_i {\mathcal N}_i
\sup_{\theta, Y}\frac{F'(|X|)e_i-F'(|X-\theta Y|)e'_i}{ F(|X|)}\,.
$$
Now we write
$$
\begin{array}{c}
F'(|X|)e_i-F'(|X-\theta Y|)e'_i=\\ [6pt]
 F'(|X|)(e_i-e_i')+(F'(|X|)-F'(|X-\theta Y|))e'_i\,.
\end{array}
$$
Since $0<\theta<1$ and $Y$ is bounded (even more, it goes to zero), and we know that
$F'(r)/F(r)\to 0$ as $r\to\infty$, this factor goes to zero uniformly in $X$ as $t\to\infty$ and the result follows. \qed


\section{Application to the Fractional Fokker-Planck Equation }\label{sec.ffpe}

Take the fractional Heat Equation posed in the whole space $\ren$ for $\tau>0$:
\begin{equation*}
u_\tau+(-\Delta_y)^s u=0\,,
\end{equation*}
with notation $u=u(y,\tau)$ that is useful since we want to save the standard notation $(x,t)$ for later use. We know the  fundamental solution
$$
P(y,\tau)=C\,\tau^{-n/2s}e^{-y^2/\tau^{1/2s}}
$$
that has a selfsimilar form. It  was  proved in previous sections that this fundamental solution is an attractor for all solutions in its {\sl basin of attraction}, consisting on all solutions with initial data that belong to $L^1(\ren)$ with integral $M=1$. See Sections  \ref{sec.conv} and \ref{conv.re}.
%

\medskip

\noindent $\bullet$ {\bf Fractional Fokker-Planck equation.} Scaling on the variables $u$ and $y$ to factor out the expected size of both which must mimic the kernel sizes, and then take logarithmic scale for the new time
\begin{equation}\label{eq.1}
 u(y,\tau)=v(x,t)\,(1+\tau)^{-n/2s}, \qquad t= \log(1+ \tau).
\end{equation}
After some simple computations this leads to the well-known \sl Fractional Fokker-Plank equation (FFPE) \rm for $v(x,t)$:
\begin{equation}\label{eq.ffp}
v_t+ (- \Delta)^{s}_x v = \frac1{2s} \nabla_x\cdot(x\,v)\,.
\end{equation}
We will keep the notation $v(x,t)$ for the solutions of this fractional Fokker-Planck equation throughout this section. We can write it as \ $v_t=L_1(v),$ where the Fokker-Planck operator $L_1=L_{FFP}$ can be written in more explicit form as
$$
L_1(v)=-(-\Delta)^{s}_x v + \frac1{2s} \nabla_x\cdot(x\,v)\,.
$$
We check now that when we look for stationary solutions by putting $v_t=0$ we
get the equation
\begin{equation}\label{eq.fundsol.1}
2s\,(-\Delta)^s F = nF+ r\partial_{r} F\,,\quad r=|x|\,,
\end{equation}
found in \eqref{eq.fundsol} for the stationary solution. As we already said, there is an explicit solution with finite mass, but no explicit formula exists for other exponents.

The theory we have established for the FHE implies immediately that the FFP equation generates a continuous contraction semigroup in all $L^p$ spaces, and also that the $L^1$ integral is conserved. The asymptotic result we are aiming at consists precisely of proving that when $v_0(x)$ is integrable with mass $M=1$ then $v(x,t)$ will tend to $F$ as $t\to\infty$. The translation of Theorem \eqref{main.convthm} and  \eqref{thm.conv.2} gives.

\begin{thm} \label{main.convthm.FFP} Let $v_0\in L^1(\ren)$ and let  $\int v_0(x)dx=M$ be its mass. Let $v(t)=v(\cdot,t)$ be the solution of the FFPE with $s\in (0,1)$, posed in the whole space. We show that
\begin{equation}\label{ffk.cr.l1}
\lim_{t\to\infty} \|v(t)-M F\|_p\to 0
\end{equation}
for all $1\le p\le \infty$. If moreover the first moment of $v_0$ is finite then
\begin{equation}\label{ffk.conv.fm}
\displaystyle  \|v(t)-M F\|\le C_1 {\mathcal N}_1(v_0)e^{-t/2s}\,,
\end{equation}
Both convergences are sharp for the intended classes of data.
\end{thm}

The results of the other theorems can also be translated. On the other hand, we know that the derivatives
$Q_i(y,\tau)=\partial_{y_i}P(y,\tau;s)$ are still solutions of equation \eqref{eq.1}. They have the form
$$
Q_i(y,\tau)=\tau^{-(n+1)/2s}G_i(y\tau^{-1/2s})
$$
A simple calculation shows that $G_i$ is an eigenfunction of operator $L_1$ with eigenvalue $1/2s$:
$L_1 G_i=-(1/2s)G_i$. This holds for all $i=1,\dots,n$. Higher derivatives allow to find higher eigenvalues.
The groundstate is $F$ with eigenvalue 0.

\medskip

\subsection{The method of entropies} The same problem has been studied by Biler and Karch in \cite{BK2003}, where they suggest the name L\'evy Fokker Planck equation, and they also consider extensions of the form
\begin{equation}
\partial_t v + {\mathcal L} v= \nabla\cdot( v\nabla V)
\end{equation}
where $\mathcal L$ is a suitable Markov diffusion operator that generates a positivity and mass preserving semigroup, and $V$ is a potential which is large enough as $|x| \to\infty$ so that it is confining, so convergence to a steady state is expected. They claim that the motivation to study the equation and its extensions stems from the probability theory where Fokker-Planck equations are deeply connected with (nonlinear) stochastic differential equations driven by Gaussian and L\'evy processes, as explained in \cite{SLDYL}.

They then use the method of entropies, nowadays in wide use (see its description in the Gaussian case in our notes \cite{Vaz17}), to get convergence with exponential decay for general integrable data, see Theorem 3.1 and Corollary 3.1 of \cite{BK2003}. This result is too good for general solutions in view of the fact that we have proved that in general there cannot be an a priori given decay rate. In fact, they make the assumption that the initial entropy is finite, which restricts the class of initial data and is compatible with our results above. The result was completed by Gentil and Imbert  \cite{GentImb} who used a related logarithmic Sobolev inequality.

Other works worth mentioning are  \cite{BKW99, BKW01, BG1960, Chafai, SLDYL}. More recently, in \cite{Trist15} asymptotic estimates with rates are obtained in weighted Lebesgue spaces $L^1(|x|^k)$. See also \cite{Toscani2016}.

\subsection{The fractional  Ornstein-Uhlenbeck semigroup}
There is a well-known connection between the classical heat equation posed in the whole Euclidean space, the Fokker-Planck equation and the Ornstein-Uhlenbeck semigroup. This is documented in the literature on Mathematical Physics, and we have reported on it in \cite{Vaz17}. In the fractional case we define the {Ornstein-Uhlenbeck operator} $L_2=L_{OU}$ by duality with respect to the Fokker-Planck operator $L_1$. For every two convenient functions $w_1$ and $w_2$ we have
\begin{equation}
\int_{\ren} (L_1 w_1)\,w_2\, dx=\int_{\ren} w_1\,(L_2w_2)\,dx.
\end{equation}
The important consequence of this computation is that $A=-L_2$ is a positive and self-adjoint operator in the Hilbert space $X=L^2(d\mu)$. This is a rather large space that includes all functions with polynomial growth. We will keep the notation $w(x,t)$ for the solutions of the Ornstein-Uhlenbeck equation throughout this section. By simple computations we get
\begin{equation}
w_t= L_2(w):=-(-\Delta)^s w - \frac1{2s} x\cdot\nabla w\,.
\end{equation}
By duality with respect to the FFP equation, we immediately conclude that the FOU equation generates a continuous contraction semigroup in all $L^p$ spaces.

The translation of our previous convergence theorems to this semigroup is an interesting question. We will not pursue further this very interesting topic here.


\section{Ideas on convergence for related nonlinear equations}\label{sec.other}

Convergence results towards selfsimilar special solutions are well known in the field of nonlinear diffusion, the main objects of study being the Porous Medium Equation and the Fast Diffusion Equation, both of the form $\partial_t u=\Delta u^m$, $m\ne 1$, and the $p$-Laplacian Equation, $\partial_t u=\nabla(|\nabla u|^{p-2}\nabla u)$, $p\ne 2$. Since the equations are nonlinear, the techniques based on linear representation in terms of a kernel cannot work.
Alternative techniques are successful and are described in the recent survey papers \cite{VazCIME} and \cite{Vaz17}.
In the last years there has been a surge of interest in the translation of the results and techniques
equations combining fractional diffusion and porous medium nonlinearities. This  gives rise to interesting mathematical models that have been studied in quite recently.

One of the models is a direct nonlinear variant of the FHE studied here. The form is
$$
\partial_t u+ (-\Delta)^s(|u|^{m-1}u)=0.
$$
It has been studied in \cite{dpqrv1,DPQRV2}.  Unique fundamental solutions of Barenblatt type were described in \cite{VazBar2014}, where convergence (without rates) was proved by a method called the scaling method. Another model was proposed by Caffarelli and the author and deals with the equation
$$
\partial_t u= \nabla(u\,\nabla(-\Delta)^{-s}).
 $$
 This model admits self-similar solutions that we may call fractional Barenblatt solutions. The entropy method is used in \cite{CV2} to establish asymptotic convergence without rates. Rates in 1D were obtained in \cite{CH2013}. Convergence with rates in several dimensions is not known.

 A very interesting extension of equation \eqref{FHE} is considered in \cite{KSZ}, where the authors consider the bifractional equation, that is, including both time fractional and space fractional derivatives. The paper contains a general study of asymptotic behaviour in that generality. The case of order 1 in time is also included the paper as a limit case and they obtain a number of statements and results overlapping to a degree with the present paper. See also the previous work \cite{VZ}.

In another direction, nonlinear fractal conservation laws were studied by a large number of authors, like \cite{Ali07, AIK, BKW99, CifJak2011, DrImb06}. They usually combine  fractal diffusion and  convection.

Quasi-geostrophic equations of the form $\partial_t u + (-\Delta)^s u+ {\bf v}\cdot\nabla u=0$, where $u$ is a scalar field and $\bf v$ a vector field related to $u$, were much studied as a possible approximation to the complete equations of viscous fluid flows, see \cite{CVass, Const2017, CoCo04}. Chemotaxis systems have been studied with nonlocal and/or nonlinear diffusion, cf. \cite{Esc06} and many later works. The asymptotic problems offer many open problems.

More equations and information are contained in the survey paper \cite{VazCIME}.

\medskip


\medskip



\noindent {\large \sc Acknowledgment}

\noindent  Funded by Project MTM2014-52240-P  (Spain).  The author would like to thank the hospitality of the Mathematical Institute of the University of Warwick. I thank R. Zacher for information about his paper \cite{KSZ} and the anonymous referee for an interesting suggestion.

\




\bibliographystyle{amsplain}

\

\

\noindent {\bf Keywords.} Fractional Heat Equation, Convergence to the fundamental solution.

\medskip

\noindent 2010 {\bf Mathematics Subject Classification}. 35B40, 35R11, 35K55.


\

\noindent e-mail address:\texttt{~juanluis.vazquez@uam.es}

\end{document}